\documentclass[finalversion]{FPSAC2017-arxiv}

\articlenumber{77}

\usepackage{pgf,tikz,kbordermatrix,xcolor}
\usetikzlibrary{arrows,shapes.geometric}
\usepackage[mathscr]{eucal}
\usepackage{color}
\usepackage{lipsum}

\newcommand{\wt}{\operatorname{wt}}
\newcommand{\conv}{\operatorname{conv}}
\newcommand{\ch}{\operatorname{ch}}

\newcommand{\g}{\mathfrak{g}}
\newcommand{\n}{\mathfrak{n}}
\newcommand{\h}{\mathfrak{h}}
\newcommand{\fl}{\mathfrak{l}}
\newcommand{\Z}{\ensuremath{\mathbb{Z}}}

\newtheorem{theo}[equation]{Theorem}

\newtheorem{cor}[equation]{Corollary}
\newtheorem{pro}[equation]{Proposition}

\crefname{theo}{Theorem}{Theorems}

\numberwithin{equation}{section}
\numberwithin{figure}{section}

\begin{document}

\title[The Weyl--Kac weight formula]{The Weyl--Kac weight formula}

\author{Gurbir Dhillon\thanks{\href{mailto:gsd@stanford.edu}{gsd@stanford.edu}. G.D.~was partially supported by the
  Department of Defense (DoD) through the NDSEG fellowship.} \and
  Apoorva Khare\thanks{\href{mailto:khare@stanford.edu}{khare@stanford.edu}}}
\address{Departments of Mathematics and Statistics,
Stanford University, Stanford, CA 94305, USA
\vspace{1ex}\begin{center}To Cora Bernard, with gratitude.\end{center}}

\received{November 14, 2016}

\revised{March 30, 2017}

\abstract{We provide the first formulae for the weights of all simple
highest weight modules over Kac--Moody algebras. For generic highest
weights, we present a formula for the weights of simple modules similar to
the Weyl--Kac character formula. For the remaining highest weights, the
formula fails in a striking way, suggesting the existence of
`multiplicity-free' Macdonald identities for affine root systems.}

\keywords{Macdonald identity, Weyl--Kac formula, Kac--Moody algebra,
highest weight module}

\maketitle

\section{Introduction}

The finite dimensional simple representations of complex
semisimple Lie algebras have long been objects of central concern in
algebraic combinatorics.  For a semisimple Lie algebra $\g$, with
triangular decomposition $\g = \n^- \oplus \h \oplus \n^+$ and Weyl group
$W$, such representations are parametrized by weights $\lambda$ in the
dominant Weyl chamber in $\h^*$. The character of the corresponding
simple module $L(\lambda)$ is given by the Weyl character formula
\cite{Weyl}. If for $\mu \in \h^*$, we write $L(\lambda)_\mu$ for the
corresponding weight space, and writing  $\Delta^+$ for the weights of
$\n^+$, i.e. the positive roots, and $\rho$ for half their sum, the
formula reads as: 
\begin{equation}\label{w}
\sum_{\mu \in \h^*} \dim L(\lambda)_\mu e^\mu = \frac{\sum_{w \in W}
(-1)^{\ell(w)} e^{w(\lambda + \rho) - \rho}}{\prod_{\alpha \in \Delta^+}
(1 - e^{-\alpha})}.
\end{equation}

We wish to remark on two points at this juncture. First, there exist
positive combinatorial formulae for $\dim L(\lambda)_\mu$, coming from
the theory of crystal bases. For example, for $\g$ of type $A$, these
characters are essentially Schur polynomials, and weight space
multiplicities are given by counting semistandard Young tableaux. Second,
for the more basic question of which $\dim L(\lambda)_\mu$ are positive,
there are two very simple descriptions. Indeed, one can (i) take the
convex hull of the Weyl group orbit $W\lambda$, and intersect this with
an appropriate translate of the root lattice to obtain the nonzero
weights. Alternatively, (ii) by $W$ invariance, it suffices to consider
dominant $\mu$, for which $L(\lambda)_\mu$ is positive if and only if
$\lambda - \mu$ is a sum of positive roots. 

When we consider a general simple highest weight module $L(\lambda)$,
i.e.~we no longer assume $\lambda$ is dominant integral, the character is
known through Kazhdan--Lusztig theory \cite{kl,bk,bb}.
For example, if $\lambda = y(\nu + \rho) - \rho$, for $\nu$ dominant
integral, and $y \in W$, we have: \begin{equation} \sum_{\mu \in \h^*}
\dim L(\lambda)_\mu e^{\mu} = \frac{ \sum_{w \in W} (-1)^{\ell(w)}
P_{wyw_\circ, yw_\circ}(1) e^{w(\lambda + \rho) - \rho}}{\prod_{\alpha
\in \Delta^+} (1 - e^{-\alpha}) }.\label{k} \end{equation} Here
$P_{x,y}(1)$ denotes the Kazhdan--Lusztig (KL) polynomial evaluated at 1,
and $w_\circ$ the longest element of $W$. What is relevant for us is that
as in the case of $y = 1$, i.e.~the Weyl character formula, there is
cancellation on the right hand side, and the situation is made more
subtle by the appearance of KL polynomials.

Unlike for finite dimensional simple modules, positive combinatorial
formulae for the multiplicities of weight spaces $\dim L(\lambda)_\mu$
are not known in general. This is a problem of longstanding interest, and
there are partial results, for example due to Mathieu and Papadopoulo in
type $A$ \cite{MP}. Surprisingly, even the simpler question of describing
which multiplicities are positive was not known until recent work of the
second author \cite{Kh1}. 

The infinite dimensional cousins of semisimple Lie algebras, the affine
Lie algebras, also play a celebrated role in algebraic combinatorics. In
broad strokes, the combinatorial representation theory of affine Lie
algebras proceeds similarly to that of semisimple Lie algebras. 

Integrable highest weight modules are the analogues of finite dimensional
simple modules, and their formal characters are given by the Weyl--Kac
formula, a modification of~\eqref{w}. For example, it was
famously realized by Kac that Macdonald's identities for affine root
systems are precisely the Weyl--Kac formula specialized to the trivial
modules of untwisted affine algebras \cite{Kac74}. As in the first remark
following~\eqref{w}, weight multiplicities in integrable highest
weight modules can be obtained positively from combinatorial descriptions
of crystal bases. A remarkable example is the Young graph of $p$-regular
partitions, which is the crystal graph of the basic representation of the
affine Lie algebra $\widehat{\mathfrak{sl}}_p$, for any prime $p$
\cite{Gr,MM}. As in the second remark following~\eqref{w}, it is
again true that the locus of weights with positive multiplicity admits
two simple descriptions similar to (i), (ii) above. 

For nonintegrable simple highest weight modules $L(\lambda)$, as in
finite type the formal characters are again largely understood. However,
the answer is considerably subtler. For integral $\lambda$ of positive or
negative level, the characters are determined by formulae similar to~\eqref{k}, now involving inverse Kazhdan--Lusztig polynomials
and Kazhdan--Lusztig polynomials for the affine Weyl group, respectively
\cite{kas,kata}. For $\lambda$ at the {\em critical} level, which plays a
distinguished role in the Geometric Langlands program
\cite{Frenkel-Gaitsgory}, the Feigin--Frenkel conjecture predicts a
formula involving {\em periodic} Kazhdan--Lusztig polynomials, and a
modified Weyl denominator \cite{afie}. As for semisimple Lie algebras, so
far individual weight multiplicities have resisted combinatorial
interpretation; moreover, even the simpler question of characterizing
which multiplicities are positive was unanswered. 

In this extended abstract of the paper \cite{DK1}, we explain a solution
to this latter problem. We will present three positive formulae for the
weights of arbitrary simple modules over semisimple and affine Lie
algebras. We will then provide a formula for the weights of all
nontrivial simple modules that is strikingly similar to the Weyl--Kac
formula, i.e.~involving signs but no Kazhdan--Lusztig polynomials or
their variants. The truth of the formula for nontrivial modules suggests
extensions of the formulae of Brion and Khovanskii--Pukhlikov for
exponential sums over lattice points of (virtual) polyhedra. The failure
of the formula for trivial modules suggests the existence of
`multiplicity-free' Macdonald identities. 

\section{Three formulae for the weights of simple modules}

Before stating our results, let us establish some context for
non-experts. Semisimple and affine Lie algebras both admit transparent
presentations by generators and relations which can be read off of the
Dynkin diagram of the corresponding root system. Such Lie algebras built
from Dynkin diagrams are known as Kac--Moody algebras. If $I$ is the set
of nodes of the Dynkin diagram, then for each $i \in I$ one has operators
$f_i, h_i, e_i$. What is relevant here is that the $h_i$ pairwise
commute, and lie in a maximal commutative Lie subalgebra $\h$. The $f_i,
e_i,$ often called lowering and raising operators, respectively, are
simultaneous eigenvectors for the  adjoint action of $\h$, with opposite
eigenvalues $-\alpha_i, \alpha_i \in \h^*$, respectively. The eigenvalues
$\alpha_i, i \in I$, are called the positive simple roots.

Let $\g$ be a Kac--Moody algebra.  A $\g$-module $V$ is called {\em
highest weight} if it can be generated by a single vector $v \in V$ that
behaves very simply under the action of `two thirds' of the generators:
(i)~$h v = (h, \lambda) v$, for some $\lambda \in \h^*$ and all $h \in
\h$,
and (ii) $e_i v = 0, \forall i \in I$.\footnote{To make condition (ii)
seem more natural, it may be helpful to note that for a finite
dimensional representation of a semisimple Lie algebra, the $e_i, i \in
I$ always act nilpotently.
Moreover, the dimension of their simultaneous kernel is the number of
simple summands of the representation. Replacing condition (ii) by a
`character' in $(\n^+/[\n^+, \n^+])^\vee$ leads to the {\em Whittaker}
modules \cite{kostant}.} Such a vector $v$ is then automatically unique
up to scalars.

As in the introduction, the simple highest weight modules over $\g$ are
parametrized by their highest weight $\lambda \in \h^*$, and we write
$L(\lambda)$ for the corresponding simple module.  We now present three
formulae for the weights of $L(\lambda)$, which are defined to be:
\begin{equation}\label{Ewts}
\wt L(\lambda) := \{ \mu \in \h^*: \dim L(\lambda)_\mu > 0 \}.
\end{equation}

The first formula uses restriction to a Levi subalgebra corresponding a
subdiagram of the Dynkin diagram. Specifically, write $I_\lambda := \{ i
\in I: (h_i, \lambda) \in \Z^{\geqslant 0} \}$. Write $\fl$ for the
subalgebra of $\g$ generated by $\h, e_i, f_i, i \in I_\lambda$. For $\nu
\in \h^*$, write $L_\fl(\nu)$ for the simple highest weight module for
$\fl$ of highest weight $\nu$. Finally, for a subset $S \subset \h^*$,
write $\Z^{\geqslant 0} S \subset \h^*$ for the set of nonnegative
integral combinations of elements of $S$. The formula then reads as:

\begin{theo}[Integrable Slice Decomposition, Dhillon--Khare
\cite{DK1}]\label{p1}
\begin{equation}
\wt L(\lambda) = \bigsqcup_{\mu \in \Z^{\geqslant 0} \{\alpha_i, i \in I
\setminus I_\lambda\} } \wt L_{\mathfrak{l}}(\lambda - \mu). \label{intsl}
\end{equation}
\end{theo}

The usefulness of \cref{p1} is that $\fl$ is the direct sum of an
abelian Lie algebra which acts semisimply on $L(\lambda)$ and the
Kac--Moody algebra associated to the Dynkin diagram with nodes
$I_\lambda$. Moreover, for the latter algebra each $L_\fl(\lambda - \mu)$
is an integrable highest weight module, for which the weights are well
known.

\begin{figure}[ht]
\begin{tikzpicture}[line cap=round,line join=round,>=triangle 45,x=1.0cm,y=1.0cm]
\coordinate (A1) at (3,6); 
\coordinate (A2) at (5,6);
\coordinate (A3) at (4,5.25);
\draw (A1) -- (A3);
\draw [fill=blue,blue] (A1) -- (A2) -- (A3) -- cycle;
\draw (5,6.6) node[anchor=north west] {$\lambda$};
\draw (3,6)-- (5,6);  
\draw (5,6)-- (4,5.25);
\draw (4,5.25)-- (3,6);
\coordinate (BH1) at (1.5,3); 
\coordinate (BH2) at (3,4);
\coordinate (BH3) at (5,4);
\coordinate (BH4) at (6.5,3);
\coordinate (BH5) at (5.5,2.25);
\coordinate (BH6) at (2.5,2.25);
\draw (BH1) -- (BH6);
\draw [fill=blue!30,blue!30] (BH1) -- (BH2) -- (BH3) -- (BH4) -- (BH5) -- (BH6) -- cycle;
\draw [dash pattern=on 6pt off 6pt] (1.5,3)-- (3,4); 
\draw [dash pattern=on 6pt off 6pt] (3,4)-- (5,4);
\draw [dash pattern=on 6pt off 6pt] (5,4)-- (6.5,3);
\draw (6.5,3)-- (5.5,2.25);
\draw (5.5,2.25)-- (2.5,2.25);
\draw (2.5,2.25)-- (1.5,3);
\coordinate (TH1) at (2.25,4.5); 
\coordinate (TH2) at (3,5);
\coordinate (TH3) at (5,5);
\coordinate (TH4) at (5.75,4.5);
\coordinate (TH5) at (4.75,3.75);
\coordinate (TH6) at (3.25,3.75);
\draw (TH1) -- (TH6);
\draw [fill=blue!55,blue!55] (TH1) -- (TH2) -- (TH3) -- (TH4) -- (TH5) -- (TH6) -- cycle;
\draw [dash pattern=on 6pt off 6pt] (2.25,4.5)-- (3,5); 
\draw [dash pattern=on 6pt off 6pt] (3,5)-- (5,5);
\draw [dash pattern=on 6pt off 6pt] (5,5)-- (5.75,4.5);
\draw (5.75,4.5)-- (4.75,3.75);
\draw (4.75,3.75)-- (3.25,3.75);
\draw (3.25,3.75)-- (2.25,4.5);
\draw (3.8,6)-- (3.95,6.15); 
\draw (3.8,6)-- (3.95,5.85);
\draw (2.8,6.7) node[anchor=north west] {$-\alpha_1 - \alpha_2$};
\draw (4.5,5.625)-- (4.52,5.85); 
\draw (4.5,5.625)-- (4.67,5.54);
\draw (4.1,5.55) node[anchor=north west] {$-\alpha_1$};
\draw [dash pattern=on 6pt off 6pt] (3,6)-- (3,1.1); 
\draw [dash pattern=on 6pt off 6pt] (5,6)-- (5,1.1);
\draw (3,6)-- (0.75,1.5); 
\draw (4,5.25)-- (1.75,0.75);
\draw (5,6)-- (7.25,1.5); 
\draw (4,5.25)-- (6.25,0.75);
\draw (6.3,3.4)-- (6.1,3.45); 
\draw (6.3,3.4)-- (6.4,3.6);
\draw (6.1,4.1) node[anchor=north west] {$-\alpha_3$};
\draw (0.2,0.5) node[anchor=north west] {$\g = \kbordermatrix{
& \alpha_1 & \alpha_2 & \alpha_3 \\
h_1 & 2 & -1 & 0\\
h_2 & -1 & 2 & -1\\
h_3 & 0 & -1 & 2}$, \quad $I_\lambda = \{ 1, 2 \}$};
\coordinate (BS1) at (11.9,5.35); 
\coordinate (BS2) at (10.5,5);
\coordinate (BS3) at (8.5,3);
\coordinate (BS4) at (10,2);
\coordinate (BS5) at (15,1);
\coordinate (BS6) at (16,1.25);
\draw (BS1) -- (BS6);
\draw [fill=blue!30,blue!30] (BS1) -- (BS2) -- (BS3) -- (BS4) -- (BS5) -- (BS6) -- cycle;
\draw [dash pattern=on 6pt off 6pt] (11.9,5.35)-- (10.5,5); 
\draw [dash pattern=on 6pt off 6pt] (10.5,5)-- (8.5,3);
\draw (8.5,3)-- (10,2);
\draw (10,2)-- (15,1);
\draw (15,1)-- (16,1.25);
\coordinate (TS1) at (12.25,6); 
\coordinate (TS2) at (10.25,5.5);
\coordinate (TS3) at (9.25,4.5);
\coordinate (TS4) at (10,4);
\coordinate (TS5) at (12.5,3.5);
\coordinate (TS6) at (14.5,4);
\draw (TS1) -- (TS6);
\draw [fill=blue,blue!55] (TS1) -- (TS2) -- (TS3) -- (TS4) -- (TS5) -- (TS6) -- cycle;
\draw [dash pattern=on 6pt off 6pt] (12.25,6)-- (10.25,5.5); 
\draw [dash pattern=on 6pt off 6pt] (10.25,5.5)-- (9.25,4.5);
\draw (9.25,4.5)-- (10,4);
\draw (10,4)-- (12.5,3.5);
\draw (12.5,3.5)-- (14.5,4);
\draw (10,6)-- (8,2); 
\draw (10,6)-- (10,1);
\draw (10,6)-- (15.5,0.5);
\draw [dash pattern=on 6pt off 6pt] (10,6)-- (12.5,1);
\draw (9.4,6.6) node[anchor=north west] {$\lambda$};
\draw (9.25,3.75)-- (9.2,3.5); 
\draw (9.25,3.75)-- (9.05,3.75);
\draw (9,3.6) node[anchor=north west] {$-\alpha_0$};
\draw (9.25,2.5)-- (9.05,2.4); 
\draw (9.25,2.5)-- (9.2,2.75);
\draw (8.5,2.4) node[anchor=north west] {$-\alpha_1$};
\draw (9.625,5.25)-- (9.8,5.35); 
\draw (9.625,5.25)-- (9.525,5.45);
\draw (8.5,5.5) node[anchor=north west] {$-\alpha_2$};
\draw (8.5,0.5) node[anchor=north west] {$\g = \kbordermatrix{
& \alpha_0 & \alpha_1 & \alpha_2 \\
h_0 & 2 & -2 & -1\\
h_1 & -2 & 2 & 0\\
h_2 & -1 & 0 & 2}$, \quad $I_\lambda = \{ 0, 1 \}$};
\end{tikzpicture}
\caption{Two illustrations of \cref{p1}, with finite and infinite
integrability. Note we only draw three of the infinitely many `integrable
slices', or rather the convex hulls of their weights, in each figure (in
the righthand figure, the top slice is a point).
Similarly, we only include four of the infinitely many rays leaving
$\lambda$ in the righthand figure.}
\label{Fig1}
\end{figure}
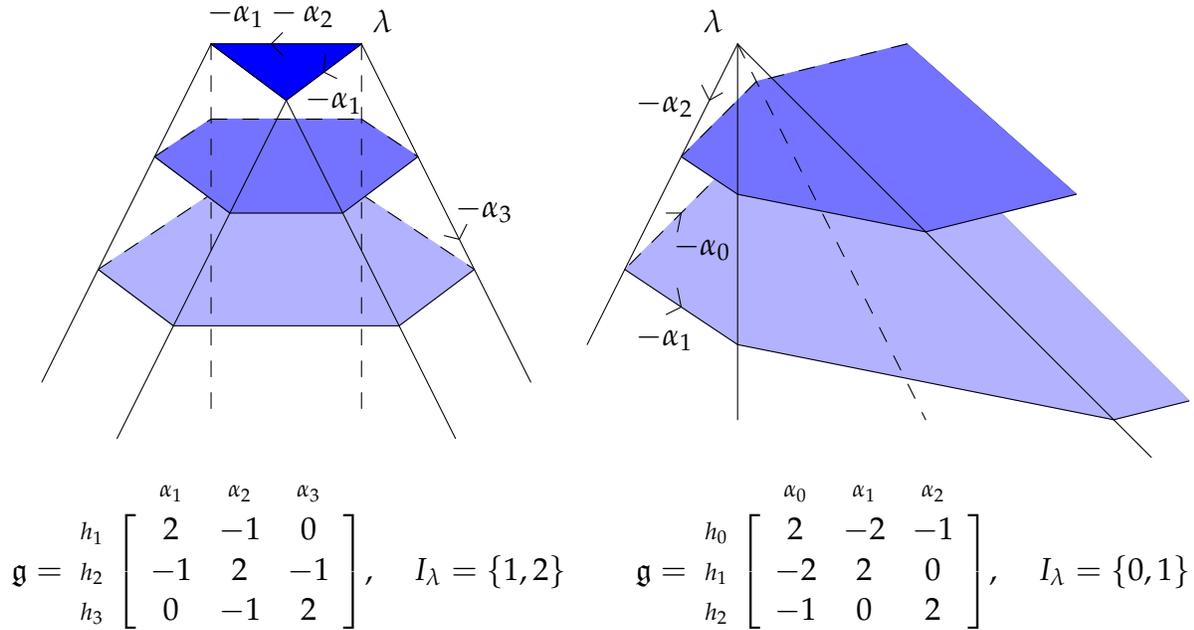

We include two illustrations of \cref{p1} in \cref{Fig1}. In
the lefthand example, the Levi $\fl$ contains a copy of $\mathfrak{sl}_3$
generated by $e_1, e_2, f_1, f_2$. As $3 \notin I_\lambda$, the weight
spaces of $L(\lambda)$ corresponding to  $\lambda, \lambda - \alpha_3,
\lambda - 2 \alpha_3,$ etc.~will all be nonzero. Moreover, they are
highest weight with respect to $\mathfrak{sl}_3$, and generate finite
dimensional $\fl$ representations. The convex hulls of weights of the
first three of these representations appear here shaded. These are the
`integrable slices' of \cref{p1}, and the theorem says these
$\mathfrak{sl}_3$ representations produce all weights of $L(\lambda)$.
The righthand example, and~\eqref{intsl} in general, can be read
similarly. 

The second formula shows the relationship between $\wt L(\lambda)$ and
its convex hull. We recall the standard partial order on weights, where
$\nu \leqslant \nu'$ if and only if $\nu' - \nu$ is a nonnegative
integral combination of positive simple roots. For a subset $X \subset
\h^*$, we write $\conv X$ for its convex hull. When $X = \wt L(\lambda)$,
we abbreviate this to $\conv L(\lambda)$.

\begin{theo}[Dhillon--Khare \cite{DK1}]\label{p3}
\begin{equation}
\wt L(\lambda) = \conv L(\lambda) \cap \{\mu \in \h^*: \mu \leqslant
\lambda \}.
\end{equation} 
\end{theo}

The question of whether the weights of simple modules $L(\lambda)$ are no
finer an invariant than their convex hull was raised by Daniel Bump.
\cref{p3} answers this question affirmatively. Note the similarity
of \cref{p3} to description (i) of the weights of simple finite
dimensional modules mentioned in the introduction.

\cref{p3} is complemented by the following description of the
convex hull of a simple highest weight module. We recall that the Weyl
group $W$ of $\g$ is generated by reflections in $\h^*$ indexed by nodes
of the Dynkin diagram $s_i, i \in I$. Let us write $W_{I_\lambda}$ for
the subgroup generated by $s_i, i \in I_\lambda$.

\begin{pro}[Dhillon--Khare \cite{DK1}] \label{ray}
\begin{equation}
\conv L(\lambda) = \conv \bigcup_{w \in W_{I_\lambda}, i \in I \setminus
I_\lambda} w\{ \lambda - \Z^{\geqslant 0} \alpha_i \}.
\end{equation} 
When $I_\lambda = I$, by the right hand side we mean $\bigcup_{w \in W}
w\lambda$. 
\end{pro}

The third formula uses the Weyl group action on $\h^*$. We will use the
following parabolic analogue of the dominant chamber:
\[
P^+_{\lambda} := \{ \nu \in \h^*: (h_i, \nu) \in \Z^{\geqslant 0},
\forall i \in I_\lambda \}.
\]
Then the final formula reads as:

\begin{theo}[Dhillon--Khare \cite{DK1}]\label{p2}
Suppose $\lambda$ has finite stabilizer in $W_{I_\lambda}$. Then:
\begin{equation}
\wt L(\lambda) = W_{I_{\lambda}} \{ \mu \in P^+_{\lambda}: \mu
\leqslant \lambda \}.
\end{equation}
\end{theo}

Note the similarity of \cref{p2} to description (ii) of the
weights of finite dimensional modules appearing in the introduction. Let
us comment on the hypothesis on $\lambda$ appearing in \cref{p2},
as it will reappear shortly.  The assumption is met by all $L(\lambda)$
for semisimple Lie algebras. For $\g$ an affine algebra with connected
Dynkin diagram, let us call a simple module $L(\lambda)$ {\em trivial} if
$\dim L(\lambda) = 1$, or equivalently $(h_i, \lambda) = 0, \forall i \in
I$. Then the assumption is met by all $L(\lambda)$ which are not trivial.  

As stated, \cref{p2} does not hold for all $\lambda$. For experts,
this can be seen by thinking about a trivial module for
$\widehat{\mathfrak{sl}}_2$ and lowering by an imaginary root. However,
as we explain in \cite{DK2}, it can be corrected by using a refinement of
the partial order $\leqslant$ introduced by Kac and Peterson
\cite{KacPet}. 

For integrable $L(\lambda)$, \cref{p1} is a tautology, and
\cref{p3,p2} are well known. \cref{p1,p3}
are due to the second named author for semisimple Lie algebras
\cite{Kh1}. All other cases are to our knowledge new. 

The above descriptions of $\wt L(\lambda)$ are particularly striking in
infinite type. For $\g$ affine, the formulae  are uniform across the
negative, critical, and positive levels, in contrast to the (partly
conjectural) character formulae discussed in the introduction. When $\g$
is symmetrizable\footnote{For non-experts, this is a technical condition
on the multiplicities of edges between nodes of the Dynkin diagram, which
is automatic for semisimple and affine Lie algebras.}, we similarly
obtain weight formulae at critical $\lambda$. The authors are unaware of
even conjectural formulae for $\ch L(\lambda)$ in this case. Finally,
when $\g$ is non-symmetrizable it is a notoriously difficult problem to
compute exact multiplicities. To wit, it is completely unknown how to
compute weight space multiplicities for integrable $L(\lambda)$ or even
the adjoint representation. Thus for $\g$ non-symmetrizable, \cref{p1,p3,p2} provide as much information on  $\ch
L(\lambda)$ as one could hope for, given existing methods.

\section{The Weyl--Kac weight formula}

We now turn to a formula for $\wt L(\lambda)$ similar to the Weyl--Kac
character formula (recall that $\wt L(\lambda)$ was defined in
\eqref{Ewts}). To orient ourselves, we first remind a slightly
nonstandard presentation of the Weyl--Kac formula in \cref{Patbo} below. This presentation is due to Atiyah and Bott for
semisimple Lie algebras \cite{atbo}. As we indicate in \cite{DK1}, it can
be straightforwardly adapted to Kac--Moody algebras.

We will need a little notation and a convention. The subalgebra $\n^+$ of
$\g$ generated by the $e_i, i \in  I$ is a semisimple $\h$-module. Write
$\Delta^+ \subset \h^*$ for the weights of $\n^+$. For $\alpha \in
\Delta^+$, write $\g_\alpha$ for the corresponding eigenspace. Next, for
$w \in W$ and $\alpha \in \Delta^+$, by $w(1 - e^{-\alpha})^{-1}$ we
mean the `highest weight' expansion, i.e.:
\begin{equation}\label{Egeom}
w \frac{1}{1 - e^{-\alpha}} := \begin{cases}
1 + e^{-w \alpha} + e^{-2 w \alpha} + \cdots, & w \alpha > 0, \\
- e^{ w\alpha} - e^{2w \alpha} - e^{3w\alpha} - \cdots, & w
\alpha < 0.
\end{cases}
\end{equation}

With these preliminaries, we may state:

\begin{pro}\label{Patbo}
Let $\g$ be a symmetrizable Kac--Moody algebra, and $L(\lambda)$ an
integrable highest weight module, i.e.~$I_\lambda = I$. Then:
\begin{equation}\label{ab}
\ch L(\lambda) = \sum_{w \in W} w \frac{e^\lambda}{\prod_{\alpha \in
\Delta^+} (1 - e^{-\alpha})^{\dim \g_\alpha}}.
\end{equation}
\end{pro}

Of course, for nonintegrable $L(\lambda)$, the character is much more
involved, as reminded in the introduction. It therefore may be surprising
that the weights of generic $L(\lambda)$ admit a similar description. To
do so, let us package the weights into a multiplicity-free character:
\[
\wt L(\lambda) = \sum_{\substack{\nu \in \h^* : \\ L(\lambda)_\nu \neq
0}} e^\nu.
\]
We then have:

\begin{theo}[Dhillon--Khare \cite{DK1}]\label{Tsimple}
Let $\g$ be an arbitrary Kac--Moody algebra.  If the stabilizer of
$\lambda$ in $W_{I_\lambda}$ is finite, then:
\begin{equation}\label{wks}
\wt L(\lambda) = \sum_{w \in W_{I_{\lambda}} } w
\frac{e^{\lambda}}{\prod_{i \in I} (1 - e^{-\alpha_i})}.
\end{equation}
\end{theo}

Note the similarity of~\eqref{ab} and~\eqref{wks}. \cref{Tsimple} was known in the case of $L(\lambda)$ integrable
\cite{Brion, Kass, Postnikov, Sc, Walton2}. All other cases are to our
knowledge new, even in finite type. 

For $\g$ a semisimple Lie algebra, the application of \cref{Tsimple} to
$L(0)$ is seen to be equivalent to the following identity, which may be
thought of as a `coordinate-free' denominator identity:

\begin{cor}[Dhillon--Khare \cite{DK1}]
Let $\Delta$ be a finite root system. Let $\Pi$ denote the set of all
bases for $\Delta$, cf.~\cite{Serre}. Then:
\begin{equation}
\prod_{\alpha \in \Delta} (1 - e^{-\alpha}) = \sum_{\pi \in \Pi}
\prod_{\beta \in \Delta \setminus \pi} (1 - e^{-\beta}).
\end{equation}
\end{cor}

\subsection{Brion's formula beyond polyhedra}

We wish to call the attention of the reader to two further problems
suggested by \cref{Tsimple}. 

Firstly, when the convex hull of $\wt L(\lambda)$ is a polyhedron,
\cref{Tsimple} may be obtained from Brion's formula. This is a
more general formula for exponential sums over lattice points of
polyhedra, due to Brion \cite{Brion} for rational polytopes and Lawrence
\cite{Lawrence} and Khovanskii--Pukhlikov \cite{KP} in general,
cf.~\cite{Barvinok-points}. We thank Michel Brion for bringing this to
our attention. For experts, one needs to observe that for regular
$\lambda$ the
associated polyhedron is Delzant. The case of singular $\lambda$ may then
be obtained via a deformation argument due to Postnikov \cite{Postnikov}.

For infinite dimensional $\g$, the convex hulls $\conv L(\lambda)$ are
rarely polyhedral, e.g.~since the Weyl group is infinite and $w\lambda$
are all extremal points. Instead, one has the following result:

\begin{theo}[Dhillon--Khare \cite{DK2}]\label{locp}
Let $\g$ correspond to a connected Dynkin diagram of infinite type, and
let $L(\lambda)$ be a nontrivial module, i.e.~$\dim L(\lambda) > 1$. Then
the following are equivalent:
\begin{enumerate}
\item The stabilizer of $\lambda$ in $W_{I_\lambda}$ is finite. 

\item  The convex hull of $\wt L(\lambda)$ is locally polyhedral, i.e.
its intersection with every polytope is a polytope. 
\end{enumerate}
\end{theo}

The combination of \cref{Tsimple,locp} suggest that
Brion's formula should be true for an appropriate class of locally
polyhedral convex sets containing both polyhedra and sets of the form
$\conv L(\lambda)$. This in particular would simplify the argument for
\cref{Tsimple} for integrable modules.

\subsection{`Multiplicity-free' Macdonald identities}

The second problem we wish to mention to the reader concerns the
extension of \cref{Tsimple} to trivial modules. For root systems
of infinite type, the trivial module $L(0)$ no longer satisfies the
condition of \cref{Tsimple}. Moreover, the stated equality no
longer always holds, and instead fails in very striking ways. For
example, we obtained the following identity by a direct calculation:

\begin{pro}[Dhillon--Khare \cite{DK1}]
Let $\g$ correspond to a Dynkin diagram of infinite type with two nodes.
Write $\Delta^+_{im}$ for the positive imaginary roots of $\g$. Then:
\begin{equation}\label{wow}
\sum_{w \in W} w \frac{1}{\prod_{i \in I} (1 - e^{-\alpha_i})} =  1 +
\sum_{\delta \in \Delta^+_{im}} e^{ -\delta}.
\end{equation}
\end{pro}

As in the Macdonald identities, i.e.~the denominator identity for affine
Lie algebras, the naive equality one would guess using only real roots
turns out to require correction terms coming from the imaginary roots.
Moreover, as we are deforming a `multiplicity-free' denominator identity
for finite root systems, i.e.~for $\wt L(\lambda)$ rather than $\ch
L(\lambda)$, the correction terms appearing in~\eqref{wow} are
insensitive to the multiplicity of the imaginary root spaces $\g_\delta$.
It would be very interesting to obtain identities similar to~\eqref{wow} for higher rank infinite root systems, i.e.~more
`multiplicity-free' Macdonald identities. 

\section{Some ingredients of the proofs}

Having explained the statements of some of the results of \cite{DK1}, we
now take a moment to highlight some new ingredients that went into their
proofs. While the above statements concern the simple highest weight
modules, they emerge from a study of general highest weight $\g$-modules
$V$.

For a fixed $\lambda \in \h^*$, one might try to classify all the modules
with highest weight $\lambda$. However, this turns out to be a daunting
task: even for semisimple Lie algebras of low rank, such as
$\mathfrak{sl}_5$, there can be infinitely many non-isomorphic highest
weight modules of highest weight $\lambda$. One therefore can try to use
invariants to distinguish between members of this profusion. The
following theorem says that several of these invariants, seemingly
different, are in fact the same:

\begin{theo}[Dhillon--Khare \cite{DK1}]\label{maintheo}
Let $\g$ be an arbitrary Kac--Moody algebra. 
Let $V$ be a highest weight $\g$-module. The following data are equivalent:
\begin{enumerate}
\item $I_V$, the integrability of $V$, i.e.~$I_V = \{i \in I: f_i \text{
acts locally nilpotently on } V \}$.

\item $\conv V$, the convex hull of the weights of $V$.

\item The stabilizer of $\ch V$ in $W$.
\end{enumerate}
\end{theo}

To our knowledge, \cref{maintheo} is new even for semisimple Lie
algebras. Before explaining how it connects to the problem of determining
the weights of simple modules, let us mention a convexity-theoretic
consequence.

\begin{cor}[Ray Decomposition, Dhillon--Khare \cite{DK1}]\label{drog}
Let $V$ be a highest weight module, and let $I_V$ be as in \cref{maintheo}. Let $W_{I_V}$ denote the subgroup of $W$ generated by
$s_i, i \in I_V$. Then:
\begin{equation}
\conv V = \conv \bigcup_{w \in W_{I_V}, i \in I \setminus I_V} w(\lambda
- \Z^{\geqslant 0} \alpha_i)
\end{equation}
When $I_V = I$, by the right hand side we mean $\bigcup_{w \in W}
w\lambda$. 
\end{cor}

\begin{figure}[ht]
\begin{tikzpicture}[line cap=round,line join=round,>=triangle 45,x=1.0cm,y=1.0cm]
\draw (3,6)-- (5,6);  
\draw (5,6)-- (4,5.25);
\draw (4,5.25)-- (3,6);
\draw (5,6.6) node[anchor=north west] {$\lambda$};
\draw (3.8,6)-- (3.95,6.15); 
\draw (3.8,6)-- (3.95,5.85);
\draw (2.8,6.7) node[anchor=north west] {$-\alpha_1 - \alpha_2$};
\draw (4.5,5.625)-- (4.52,5.85); 
\draw (4.5,5.625)-- (4.67,5.54);
\draw (4.1,5.55) node[anchor=north west] {$-\alpha_1$};
\draw [->,blue] [dash pattern=on 6pt off 6pt] (3,6)-- (3,1.1); 
\draw [->,blue] [dash pattern=on 6pt off 6pt] (5,6)-- (5,1.1);
\draw [->,blue] (3,6)-- (0.75,1.5); 
\draw [->,blue] (4,5.25)-- (1.75,0.75);
\draw [->,blue] (5,6)-- (7.25,1.5); 
\draw [->,blue] (4,5.25)-- (6.25,0.75);
\draw (6.3,3.4)-- (6.1,3.45); 
\draw (6.3,3.4)-- (6.4,3.6);
\draw (6.1,4.1) node[anchor=north west] {$-\alpha_3$};
\draw [dash pattern=on 6pt off 6pt] (2.25,4.5)-- (3,5); 
\draw [dash pattern=on 6pt off 6pt] (3,5)-- (5,5);
\draw [dash pattern=on 6pt off 6pt] (5,5)-- (5.75,4.5);
\draw (5.75,4.5)-- (4.75,3.75);
\draw (4.75,3.75)-- (3.25,3.75);
\draw (3.25,3.75)-- (2.25,4.5);
\draw [dash pattern=on 6pt off 6pt] (1.5,3)-- (3,4); 
\draw [dash pattern=on 6pt off 6pt] (3,4)-- (5,4);
\draw [dash pattern=on 6pt off 6pt] (5,4)-- (6.5,3);
\draw (6.5,3)-- (5.5,2.25);
\draw (5.5,2.25)-- (2.5,2.225);
\draw (2.5,2.225)-- (1.5,3);
\draw [->,blue] (10,6)-- (8,2); 
\draw [->,blue] (10,6)-- (10,1);
\draw [->,blue] (10,6)-- (15.5,0.5);
\draw [->,blue] [dash pattern=on 6pt off 6pt] (10,6)-- (12.5,1);
\draw (9.4,6.6) node[anchor=north west] {$\lambda$};
\draw (9.25,3.75)-- (9.2,3.5); 
\draw (9.25,3.75)-- (9.05,3.75);
\draw (9,3.6) node[anchor=north west] {$-\alpha_0$};
\draw (9.25,2.5)-- (9.05,2.4); 
\draw (9.25,2.5)-- (9.2,2.75);
\draw (8.5,2.4) node[anchor=north west] {$-\alpha_1$};
\draw (9.625,5.25)-- (9.8,5.35); 
\draw (9.625,5.25)-- (9.525,5.45);
\draw (8.5,5.5) node[anchor=north west] {$-\alpha_2$};
\draw [dash pattern=on 6pt off 6pt] (12.25,6)-- (10.25,5.5); 
\draw [dash pattern=on 6pt off 6pt] (10.25,5.5)-- (9.25,4.5);
\draw (9.25,4.5)-- (10,4);
\draw (10,4)-- (12.5,3.5);
\draw (12.5,3.5)-- (14.5,4);
\draw [dash pattern=on 6pt off 6pt] (11.9,5.35)-- (10.5,5); 
\draw [dash pattern=on 6pt off 6pt] (10.5,5)-- (8.5,3);
\draw (8.5,3)-- (10,2);
\draw (10,2)-- (15,1);
\draw (15,1)-- (16,1.25);
\end{tikzpicture}
\caption{Ray Decomposition, with finite and infinite integrability; here
$\g$ and $I_V$ are as in \cref{Fig1}.}
\label{Fig4}
\end{figure}
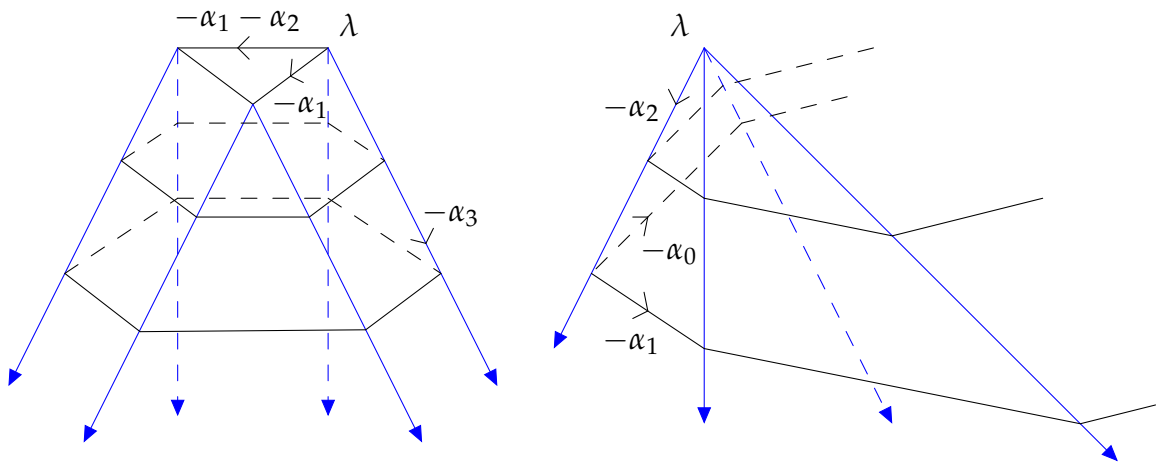

We include two illustrations of \cref{drog} in \cref{Fig4}. To our knowledge \cref{drog} was not previously
known for non-integrable modules in either finite or infinite type.

Notice that \cref{ray} is a special case of \cref{drog}. The above presentation of the convex hull can be understood
as follows. Consideration of the nodes of $I \setminus I_V$, and the
representation theory of $\mathfrak{sl}_2$ tell us that $\lambda -
\Z^{\geqslant 0} \alpha_i$ lies in $\ch V$, $\forall i \in I \setminus
I_V$. Consideration of the nodes of $I_V$ tell us the character, whence
the convex hull, should be $W_{I_V}$ invariant. The content of \cref{drog} is that these two {\em a priori} estimates are in fact enough
to generate the convex hull. 

As \cref{drog} indicates, \cref{maintheo} has useful
implications in the convexity theoretic study of highest weight modules.
In fact, in the companion work \cite{DK2}, we determine the face posets
of all convex hulls of highest weight modules. To our knowledge, this
classification had not been fully achieved even for semisimple Lie
algebras.  

Finally, let us explain how to use \cref{maintheo} to obtain
information about $\wt L(\lambda)$, i.e.~how to pass from convex hulls to
weights. For a highest weight module $V$, define the {\em potential
integrability} of $V$ to be $I^p_V := I_{\lambda} \setminus I_V$. To
justify the terminology, note these are precisely the simple directions
whose actions become integrable in quotients of $V$.

In the following theorem, we use the {\em parabolic Verma modules}
$M(\lambda, J)$, for $J \subset I_\lambda$. These are characterized by
the property that they have integrability $J$ and admit a morphism to
every highest weight module $V$ of highest weight $\lambda$ with $I_V =
J$.

\begin{theo}[Dhillon--Khare \cite{DK1}]\label{wts}
Let $V$ be a highest weight module, $V_\lambda$ its highest weight line.
Let $\fl$ denote the Levi subalgebra generated by $\h, e_i, f_i, i \in
I^p_V$. Then $\wt V = \wt M(\lambda, I_V)$ if and only if $\wt
U(\mathfrak{l}) V_\lambda = \wt U(\fl) M(\lambda, I_V)_\lambda$,
i.e.~$\wt U(\fl) V_\lambda = \lambda - \Z^{\geqslant 0} \{\alpha_i, i \in
I^p_V\}$.

In particular, if $V = L(\lambda)$, or more generally $|I^p_V| \leqslant
1$, then $\wt V = \wt M(\lambda, I_V)$.
\end{theo}

\cref{wts} is new in both finite and infinite type. \cref{maintheo,wts}, combined with some analysis of the parabolic
Verma modules, yield proofs of the results advertised in the second and
third sections. In fact, note that \cref{p1,p3,p2} correspond to the manifestations 1., 2., and 3., respectively,
of the invariant studied in \cref{maintheo}.
 
\acknowledgements{We would like to thank Galyna Dobrovolska
and the referees for their careful reviews of an earlier draft of this
manuscript, as well as their suggestions on improving the exposition. }




%
%




\end{document}